\documentclass[smallextended]{svjour3}       

\usepackage{graphicx}
\usepackage{caption}
\usepackage{amsmath,amsfonts,amssymb}
\usepackage{natbib}
\usepackage{color}
\usepackage[colorlinks=true,linkcolor=darkblue,citecolor=darkblue,
            filecolor=darkblue,urlcolor=darkgreen]{hyperref}

\definecolor{darkred}{rgb}{0.5,0.1,0.1}
\definecolor{darkgreen}{rgb}{0.1,0.5,0.1}
\definecolor{darkblue}{rgb}{0.2,0.2,1.0}

\newcommand{\eql}{\begin{equation}\label}
\newcommand{\eqn}[1]{(\ref{#1})}
\newcommand{\goto}{\rightarrow}

\newcommand{\Fig}[1]{Figure~\ref{fig:#1}}

\newenvironment{choice}{\left\{ \begin{array}{ll}}{\end{array}\right.}
\def\choose{\begin{choice}}
\def\when{&\mbox{if}}

\def\xiLAT{\xi_{\scriptscriptstyle LAT}}
\def\xiHAT{\xi_{\scriptscriptstyle HAT}}
\def\xiMSL{\xi_{\scriptscriptstyle MSL}}
\def\xiMHW{\xi_{\scriptscriptstyle MHW}}
\def\xiMHHW{\xi_{\scriptscriptstyle MHHW}}
\def\xiMLW{\xi_{\scriptscriptstyle MLW}}
\def\xiMLLW{\xi_{\scriptscriptstyle MLLW}}
\def\erf{\text{erf}}

\begin{document}
\title{The Pattern Method for Incorporating Tidal Uncertainty Into
Probabilistic Tsunami Hazard Assessment (PTHA)\thanks{This
work was done as a Pilot Study funded by BakerAECOM. Partial funding
was also provided by NSF grant numbers DMS-0914942 and DMS-1216732 of the second author.}}

\titlerunning{The Pattern Method for Tidal Uncertainty}

\author{Loyce M. Adams \and Randall J. LeVeque \and Frank I. Gonz\'alez}

\authorrunning{Adams, LeVeque, and Gonz\'alez}
\institute{University of Washington \at
              Dept. Applied Mathematics \\
              \email{lma3@uw.edu}           
           \and
           University of Washington \at
              Dept. Applied Mathematics \\
              \email{rjl@uw.edu}           
           \and
           University of Washington \at
              Dept. Earth \& Space Sciences \\
              \email{figonzal@uw.edu}           
}
\date{March 4, 2014}
\maketitle
\begin{abstract}
In this paper we describe a general framework for 
incorporating tidal uncertainty into probabilistic tsunami hazard assessment
and propose the Pattern Method and a simpler special
case called the $\Delta t$ Method as effective approaches.  
The general framework also covers the method developed by \cite{Mofjeld2007}
that was used for the 2009 Seaside, Oregon probabilistic
study by \cite{Gonzalez:2009p452}.
We show the Pattern Method is superior to past approaches
because it takes advantage of our ability to run the tsunami simulation at
multiple tide stages and uses the time history of flow depth at strategic
gauge locations to infer the temporal pattern of waves that is unique to
each tsunami source.  Combining these patterns with knowledge
of the tide cycle at a particular location improves the ability
to estimate the probability that a wave will arrive at a time when the tidal
stage is sufficiently large that a quantity of interest such as the maximum
flow depth exceeds a specified level.
\end{abstract}
\section{Introduction}\label{sec:tideintro}

A numerical tsunami modeling code typically takes as input the
seafloor deformation due to an event (such as an earthquake or
submarine landslide) and then simulates the resulting tsunami
generation and propagation.  Often the desired output is a map of
the maximum depth of flooding, maximum flow velocity, or some other quantity
of interest over some region along the coast.  We will use $\zeta(x,y)$ to
denote some generic Quantity of Interest (QoI) that might be computed as a
function of spatial location (typically maximizing some non-negative quantity
over the entire duration of the tsunami event).
We use the open source GeoClaw model (\url{http://www.clawpack.org/geoclaw})
described in detail by \cite{BergerGeorgeLeVequeMandli}
and \cite{LeVequeGeorgeBerger:an11},
but the methodology developed here could be used in conjunction with 
other tsunami modeling codes.

One use of such models is in performing Probabilistic Tsunami Hazard
Assessment (PTHA), in which a probability distribution is specified
over a space of possible seafloor deformations and the desired
output is a probabilistic map of the QoI:  At each point $(x,y)$ in the
region of interest we wish to estimate the annual probability that $\zeta(x,y)$
will exceed some value $\hat\zeta$, typically for a specified
set of {\em exceedance values} 
$\hat\zeta = \zeta_1,~\zeta_2,~\ldots,~\zeta_I$.
A curve that shows the annual probability of exceedance $P[\zeta(x,y) >
\hat\zeta]$ as a function of $\hat\zeta$ is called a {\em hazard curve}, and
plotting contours derived from these hazard curves constructed at each point
on a dense grid of $(x,y)$ locations gives the desired hazard map.

There are many uncertainties that must be taken into account in performing
PTHA.  The largest source of epistemic uncertainty is the paucity of
knowledge of the proper probability distribution for seafloor deformations.
PTHA associated with subduction zone megathrust events (e.g. Sumatra 2004,
Maule 2010, Tohoku 2011) is typically done by first developing a finite set of ``characteristic
earthquakes'' that are thought to be representative of possible events and
assigning an annual probability of occurrence to each.  If this 
correctly described the probability distribution, and if there were no other
sources of uncertainty, then computing the QoI via a single tsunami simulation 
for each event would give a set of values that could be easily combined into the
desired hazard curves.  

There are many other sources of uncertainty that must also be considered in a
full PTHA.  In this paper we consider one important source of aleatoric
uncertainty: the fact that for hypothetical tsunamis in the future we do not
know whether the tsunami will reach the region of interest at high tide or
low tide and the effect can be considerably different in some cases depending 
on the tide stage.  We assume the region of interest is sufficiently small
(e.g. a single harbor or city) that the tide can be described by a single
function $\xi(t)$ that is measured in meters relative to Mean Sea Level (MSL)
in this region.  We also assume we know $\xi(t)$ for all $t$ spanning a
sufficiently long time period (e.g. one year) that we can assume this time
series properly describes the statistical properties of the tide.  
Tide gauge records are available at many locations that can be used, 
or $\xi(t)$ can be determined analytically from Fourier series with known
coefficents for the tidal constituents,
which have also been determined for many locations.  With this assumption,
there is no epistemic uncertainty in the tide and we have only the aleatoric
uncertainty associated with the fact that the earthquake could happen at any
time, which means that the time $t_0$ of the first arrival of a tsunami at
the region of interest can be viewed as a random variable that is uniformly
distributed over the time period for which we have the representative tidal
record $\xi(t)$.

The tide function $\xi(t)$ represents
the tide at the location of
interest as a function of time, relative to Mean Sea Level (MSL) which is
taken to be $\xiMSL = 0$.  
Certain site-specific values will be referred to below and we summarize
these here, adapted from
\url{http://tidesandcurrents.noaa.gov/datum_options.html}. 
Typically the tide function $\xi(t)$ exhibits two high tides each
tidal day, one of which is often considerably higher than the other. The
value  $\xiMHW$ (Mean High Water)
denotes the average of all the high water heights observed over the National
Tidal Datum Epoch, while $\xiMHHW$ (Mean Higher High Water)
denotes the  average of the higher high water height of each tidal day.  
Similarly, $\xiMLW$ and $\xiMLLW$ denote Mean Low Water and Mean Lower Low
Water, respectively.  Finally $\xiLAT$ and $\xiHAT$ denote the lowest and
highest predicted astronomical tide expected to occur at the site. (Observed
tides may be higher or lower due to non-astronomical effects such as storm
surge or other meterological effects, for example).  
Values for Crescent City are given in Section \ref{sec:CCtides}.

In this paper we focus entirely on the following question: Given that a
particular tsunami-generating event $E$ occurs, what is the probability that
$\zeta(x,y)$ exceeds some specified level $\hat\zeta$?   We are not
concerned with estimating the probability that $E$ occurs (or any
sources of uncertainty other than the fact that $t_0$ is a random
variable distributed uniformly as described above), and so we are
really concerned with the conditional probability $P[\zeta>\hat\zeta \,|\,E]$,
where the only randomness is in
$t_0$.  Note that if $P_E$ is the probability that event $E$ occurs,
then the annual probability of exceeding $\hat\zeta$ due to this
event is then given by the product $P_E P[\zeta>\hat\zeta \,|\, E]$ and it is
these probabilities that are combined with those from other events
to obtain the hazard curves.  Henceforth we simplify notation by
simply writing $P[\zeta > \hat\zeta]$ for the conditional probability
and assume we are focusing on one possible event.  Note that we
also drop the explicit dependence on $(x,y)$ for brevity, but these
conditional probabilities will vary in space and must be computed
separately at each point of interest.

\section{Background and context}\label{sec:backcontext}
The methods and theory presented
in this paper were first developed in conjunction with
a PTHA study \citep{GonLevAda} of Crescent City, California
funded by BakerAECOM and motivated by FEMA's desire to improve products of
the FEMA Risk Mapping, Assessment, and Planning (Risk MAP) Program.
At Crescent City, the difference in tide level between Mean Lower
Low water (MLLW) and Mean Higher High water (MHHW) is about 2.1 meters.
Coastal sites with such a significant tidal range experience
tsunami/tide interactions that are an important factor in the degree
of flooding.  For example, \cite{Kowalik} conducted
a modeling study that focused on two sites, Anchorage and Anchor
Point, in Cook Inlet, Alaska. They found tsunami/tide interactions
to be very site-specific, with strong dependence on local bathymetry
and coastal geometry, and concluded that the tide-induced change
in water depth was the major factor in tsunami/tide interactions.
Similarly, a study of the 1964 Prince William Sound tsunami
\citep{Zhang} compared simulations conducted with and without
tide/tsunami interactions. They also found large site-specific
differences and determined that tsunami/tide interactions can account
for as much as 50\% of the run-up and up to 100\% of the inundation.
Thus, probabilistic tsunami hazard assessment (PTHA) studies must
account for the uncertainty in tidal stage during a tsunami event.

\cite{Houston} developed probabilistic tsunami inundation
predictions that included tidal uncertainty for points along the
US West Coast. The study was conducted for the Federal Insurance
Agency, which needed such assessments to set federal flood insurance
rates.  They considered only far-field sources in the Alaska-Aleutian
and Peru-Chile Subduction Zones, because local West Coast sources
such as the CSZ (Cascadia Subduction Zone) and Southern California Bight
landslides,
\url{http://www.usc.edu/dept/tsunamis/2005/pdf/GRL2004BorreroEtal.pdf}, 
had not
yet been discovered, and assigned probabilities to each source
based on the work of \cite{Soloviev}.  Maximum runup estimates were
made at 105 coastal sites rather than from actual inundation
computations on land.  The tidal uncertainty methodology began with
a modeled 2-hour tsunami time series that was extended by 24-hours by
appending a sinusoidal wave with an amplitude that was 40\% of the
maximum modeled wave, to approximate the observed decay of West
Coast tsunamis.  This 24-hour tsunami time series was then added
sequentially to 35,040 24-hour segments of a year-long record of
the predicted tides, each segment being temporally displaced by 15
minutes.  Determination of the maximum value in each 24-hour segment
then yielded a year-long record of maximum combined tide and tsunami
elevations, each associated with the probability assigned to the
corresponding far-field source.  Ordering the elevations and,
starting with the largest elevations, summing elevations and
probabilities to the desired levels of 0.01 and 0.002, produced the
100-year and 500-year elevations, respectively.

\cite{Mofjeld2007}
developed a tidal uncertainty methodology that, unlike that of
\citep{Houston}, does not use modeled tsunami time series.  Instead,
a family of synthetic tsunami series are constructed, each with a
period in the tsunami mid-range of 20 minutes and an initial amplitude
ranging from 0.5 to 9.0 m that decreases exponentially with the
decay time of 2.0 days, as estimated by \cite{VanDorn} for
Pacific-wide tsunamis.  As in \citep{Houston}, linear
superposition of tsunami and tide is assumed and the time series
are added sequentially to a year-long record of predicted tides at
progressively later arrival times, in 15 minute increments.  Direct
computations are then made of the probability density function (PDF)
of the maximum values of tsunami plus tide.  The results are then
approximated by a least squares fit Gaussian expression that is a
function of known tidal constants for the area and the computed
tsunami maximum amplitude; for this reason, we refer to this approach
as the Gaussian method, or G Method.
This expression provides a convenient means of
estimating the tidal uncertainty, and was used by
\cite{Gonzalez:2009p452} in the PTHA study of Seaside, OR.

In this paper we present a unified framework that will be seen to include the 
G Method used by \cite{Mofjeld2007} and \cite{Gonzalez:2009p452}.  
We also present the Pattern Method which falls within this unified framework
but has the following improvements on that methodology:
(a) The assumption of linear superposition of the tide
and tsunami waves is replaced by a methodology that utilizes multiple
runs at different tidal stages; thereby introducing nonlinearities
in the inundation process that are not accounted for in previous
methods. (b) Synthetic time series are replaced by the actual
time series computed by the inundation model.  
(c) The Pattern Method takes account of temporal wave patterns that are
unique to each tsunami source; for example, some sources produce
only one large wave, others a sequence of equally dangerous waves
that arrive over several hours.  Combining these patterns with
knowledge of the tide cycle at a particular location like Crescent
City improves estimates of the probability that a wave will arrive
at a time when the tidal stage is sufficiently large that inundation
above a level of interest occurs.  
A special case of the Pattern Method that we call the $\Delta t$ Method is
discussed first since it is easier to understand and
implement, and may be sufficient for many tsunami studies.

In Section \ref{sec:overview}, we give an overview of our framework for
calculating $P(\zeta > \hat\zeta)$ assuming event $E$ has occurred.
The $\Delta t$ Method, the Pattern Method, and the G Method are introduced
and unified under this framework in Section \ref{sec:methods}.  In Section
\ref{sec:tides_compare}, we use results from an example Crescent City PTHA
study to compare these methods.

\section{Overview of the framework}\label{sec:overview}

Recall that we consider one
specific (hypothetical) tsunami event and one location $(x,y)$ and
are attempting to compute the
conditional probability that the QoI $\zeta(x,y)$
will exceed some value $\hat \zeta$, given that this event occurs.  
In practice we estimate these only for a set of discrete {\em exceedance
values} $\hat\zeta = \zeta_i$ for $i=1,~2,~\ldots, I$, but the methodology
can be described as a general approach to determining 
a complementary cumulative distribution function
(CCDF) $\Psi(\hat\zeta)$ such that 
\eql{Psi}
P[\zeta > \hat\zeta] \approx \Psi(\hat\zeta).
\end{equation} 
Recall that the random variable is the time $t_0$ at which the tsunami arrives,
which is assumed to be uniformly distributed over a typical year of the tidal
record, for example. All statements about probability are with respect to
this underlying uniform distribution.

Once a specific tsunami event has been specified and a numerical method
chosen to estimate the QoI, the only free parameter is the (static) tide
stage $\hat \xi$ used to run the code.  We assume it is possible to run the
code for any choice of $\hat \xi$ and so in principle there is a function
$\hat Z$ so that $\hat\zeta = \hat Z(\hat\xi)$ is the value of the QoI that results from
running the code with tide stage $\hat\xi$.  In practice we cannot determine
this function for all $\hat\xi$ in finite time, but we can approximate the
function by various approaches and we refer to the approximation as
$Z(\hat\xi)$.  In this paper we consider two possibilities:
\begin{itemize}
\item Run the code at a single tide level, for example taking a nearly
worst-case value $\hat\xi=\xiMHHW$, and then assume that for other values of
$\hat\xi$ the value of $Z(\hat\xi)$ varies in some specified manner.  If the
QoI is maximum depth of inundation, then assuming linear variation with slope
1 might be the natural choice. With this choice, any change in tide level is
simply added to the inundation and 
\eql{Zlinear}
Z(\hat\xi) = \hat Z(\xiMHHW) + (\hat\xi - \xiMHHW).
\end{equation} 
This may be the only option when using a numerical model that does not
allow adjusting the sea level parameter, and was used in the 
methodology developed in \citep{Mofjeld2007}, as discussed further in
Section \ref{sec:Gmethod}.
\item Run the code at several values of $\hat\xi$ and use piecewise linear
interpolation to approximate $Z(\hat\xi)$ for intermediate values.
If the set of $\hat\xi$ values used to approximate $Z$ does not span the full
range of possibilities from $\xiLAT$ to $\xiHAT$, then it may still be
necessary to use linear extrapolation beyond the largest $\hat\xi$ used, for
example.
\end{itemize} 

\begin{figure}[t]
\hfil\includegraphics[width=2.3in]{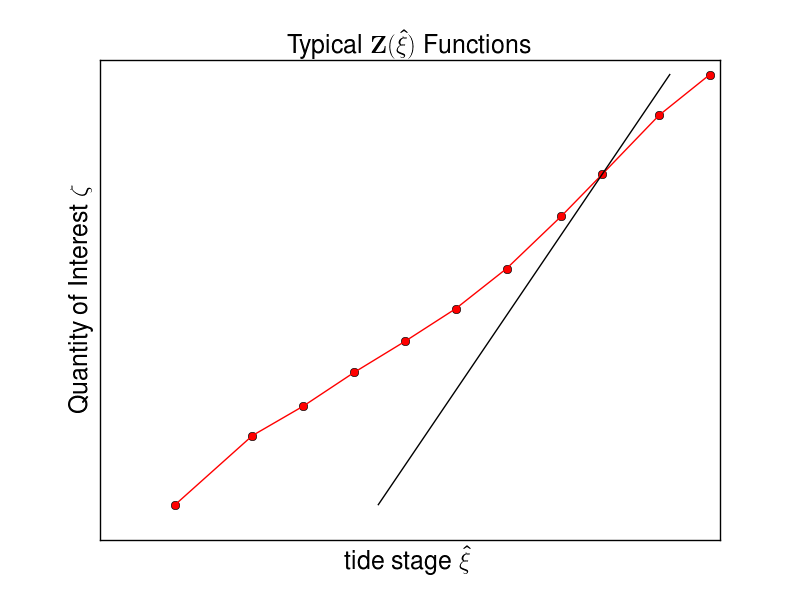}\hfil
\hfil\includegraphics[width=2.3in]{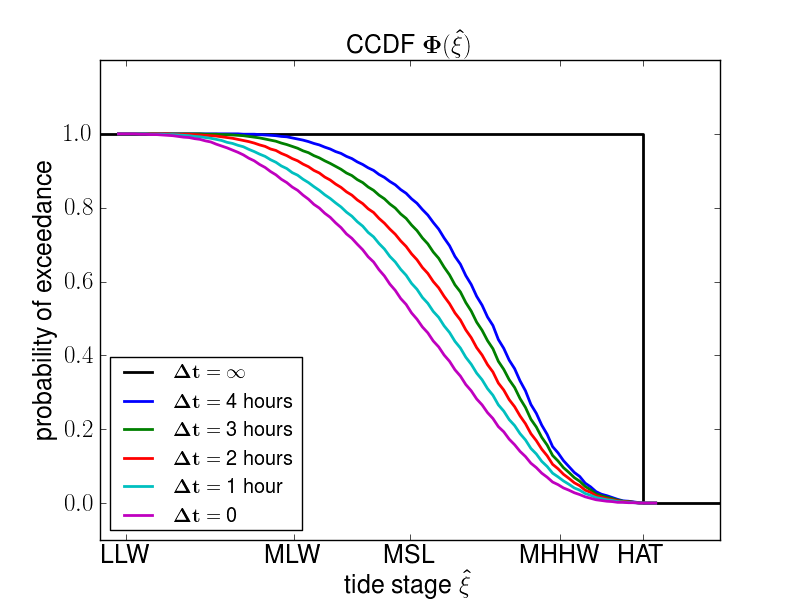}\hfil
\vskip -3mm
\caption{\label{fig:PhiDt} Left: Two possible functions $Z (\hat\xi)$, a piecewise linear
approximation (red) and a linear approximation based on a single run (black),
Right: $\Phi_{\Delta t} (\hat\xi)$} 
\end{figure}

The second approach is preferable when possible, and we have found in our study of
Crescent City that the relation between $\hat\xi$ and the maximum inundation
can be very different from \eqn{Zlinear} in many onshore regions.
Typical $Z$ functions for these two approaches are illustrated in the left plot
in Figure~\ref{fig:PhiDt}.

Our methodology also requires the inverse function $Z^{-1}$, with the
interpretation that $\hat\xi = Z^{-1}(\hat\zeta)$ is the minimal tide level
above which $\zeta$ exceeds $\hat\zeta$.  
If the function $Z(\hat\xi)$ is monotonically increasing then $Z^{-1}$ is
truly the inverse function.  If $Z$ is non-monotone, then several tide
levels $\hat\xi$ could result in the same value of $\zeta$. As a conservative
choice we might want to choose the smallest tide level above which
the QoI exceeds level $\hat\zeta$ (i.e. the infimum of the open set of
points where $Z(\hat\xi) > \hat\zeta$), and so we could define
\eql{Zinv}
Z^{-1}(\hat\zeta) \equiv \inf(\hat\xi: Z(\hat\xi) > \hat\zeta)
\end{equation} 
to make this function well-defined.  (For some quantities of interest
such as the maximum fluid velocity we have found that $Z(\hat\zeta)$
may be far from monotone and there are also other alternatives, discussed
briefly at the end of this section).

In addition to the functions $Z$ and $Z^{-1}$, the second main component of
our general methodology is another CCDF,
$\Phi(\hat\xi)$, that allows us to map a specific tide level $\hat\xi$ to the
probability that the tide will be above this level when the tsunami occurs.
There are many approaches to defining this function, depending on how one
interprets the phrase ``when the tsunami occurs''.  To explain the basic
idea we first consider the simplest approximation, which is to assume that the
tsunami consists of a single destructive wave that inundates and retreats
over a much faster time scale than the rise and fall of the tide.  In this
case, we could choose $\Phi$ to be the function 
\eql{Phi0}
\Phi_0(\hat\xi) = P[\xi(t_0) > \hat\xi]
\end{equation}
where $t_0$ is a random time, sampled for example from a year of tidal
records at the location of interest with a uniform distribution. This
function $\Phi_0$ is easily approximated from available tide gauge records at
many locations, or could be computed from $\xi(t)$ specified analytically
from the tidal constituents, which have also been determined for many
locations.  

It is also sometimes useful to discuss the probability density 
$\phi(\hat\xi) = -\Phi'(\hat\xi)$.  Note that for the CCDF $\Phi_0$ of
\eqn{Phi0}, the corresponding density $\phi_0$ has the property that
\[
\int_{\hat\xi_1}^{\hat\xi_2} \phi_0(\hat\xi)\, d\hat\xi = P[\hat\xi_1 <
\xi(t_0) < \hat\xi_2],
\]
where again $t_0$ is the  uniformly distributed
random variable and $\xi(t)$ is the known tidal variation function for the
region of interest.

We can now explain how $\Psi$ in \eqn{Psi} is determined by our approach.
If we accept $\Phi_0(\hat\xi)$ as giving the 
probability that the tide will be above level $\hat\xi$ when a very short
duration tsunami arrives, and if the tide level must be above some value
$Z^{-1}(\hat\zeta)$ in order for the QoI to exceed $\hat\zeta$, then clearly
the probability of exceedance is $\Phi_0\left(Z^{-1}(\hat\zeta)\right)$.
This gives the definition of $\Psi(\hat\zeta)$ we desire for \eqn{Psi}.
This is the key idea of our methodology.  Several variants will be
discussed in more detail based on our desire to improve on the choice $\Phi_0$.

Typically a tsunami does not consist of a single wave over a very short time
duration, but rather a series of waves arriving over the course of
several hours, during which the tide varies.  In this case accurate modeling
of the QoI might require a numerical model that also models the rise and fall
of the tide and the resulting tidal currents and how they interact with the
tsunami.  Some work has been done in this direction, see \citep{androsov_tsunami_2011}
and \citep{kowalik_tide-tsunami_2006},
but tsunami models
currently in use for forcasting or hazard assessment work do not have this
capability, even for modeling a specific event when the arrival time $t_0$ is
known.  Performing a probabilistic study where $t_0$ is random would be even
more difficult since the model would have to be run with many different
choices of $t_0$ to explore the full range of possibilities.  Instead we
focus on ways to improve the analysis in the practical case where the model
can be run at different static tide levels $\hat\xi$ but not with dynamically
varying tides.

We can still improve on the choice $\Phi_0$ in a number of ways, several of
which are explored in this paper.  If we know that destructive waves arrive
over a time period of length $\Delta t$ then it is natural to consider the
probability that the maximum value of
$\xi(t)$ is above some specified level $\hat\xi$ over a random time interval
$t_0 \leq t \leq t_0+\Delta t$, which suggests the CCDF
\eql{PhiDt}
\Phi_{\Delta t}(\hat\xi) = P\left[\max_{t_0 \leq t \leq t_0+\Delta t} \xi(t) >
\hat\xi\right].
\end{equation} 
Again $t_0$ is the random variable, uniformly distributed over 1 year, say.
Note that for $\Delta t = 0$ this reduces to $\Phi_0$ defined in \eqn{Phi0}.

Several $\Phi_{\Delta t}$ curves are illustrated in Figure~\ref{fig:PhiDt}.
The lowermost curve is $\Phi_0$ and as we increase the length of the time
interval, the probability that $\xi(t)$ will be above a fixed $\hat\xi$ value
over a random interval of this length will increase.  The limiting curve as
$\Delta t \goto \infty$ is the discontinuous piecewise constant function
\eql{PhiInf}
\Phi_\infty(\hat\xi) = 
\choose 1 \when~~ \hat\xi < \xiHAT, \\
        0 \when~~ \hat\xi \geq \xiHAT. \end{choice}
\end{equation} 
This results from the fact that over a sufficiently
long time interval we are almost surely going to observe $\xi(t)$ above any
value, up to the highest possible value that can be observed.

Choosing one of these $\Phi_{\Delta t}$ for $\Phi$ gives the $\Delta t$ Method
described further in Section \ref{sec:dtmethod}.  These CCDF's are also easily computed from
tide gauge data or tidal constituents and for many tsunami events this may be a
reasonable approach, estimating $\Delta t$ by examining the wave pattern
observed in simulations.  

For some events, however, there may be many destructive waves that arrive 
over the course of many hours or even several days, but interspersed by
periods of no waves or outflow.  
In this case taking $\Delta t$ sufficiently large to
capture all the waves might overestimate the probability that waves will
arrive when the tide is high.
So rather than looking at a single time
interval of length $\Delta t$ starting at a random time $t_0$, 
it may be more accurate to specify a {\em pattern} of disjoint time intervals
during which the large waves arrive, and then slide this pattern over the full
tidal record to determine the probability that $\max_{t} \xi(t)$ 
will be above $\hat\xi$ (for $t$ ranging over this disjoint
collection of time intervals, starting at some random time $t_0$).
This approach can be made more general by also incorporating information
about the relative magnitude of waves in the different intervals.  
This Pattern Method is discussed in more detail in Section \ref{sec:tides_pattern}.

Yet another approach to choosing $\Phi(\hat\xi)$ would be to observe that
most of the curves in \Fig{PhiDt} resemble plots of the complementary error
function (the CCDF corresponding to a Gaussian PDF) and choose
\eql{PhiErf}
\Phi(\hat\xi) = \frac 1 2 \left(1 - \erf\left((\hat\xi -
\xi_0)/\sqrt{2}\sigma\right)\right),
\end{equation} 
for some choice of the mean $\xi_0$ and standard deviation $\sigma$.  
In Section \ref{sec:Gmethod} we show that the method proposed in \citep{Mofjeld2007}
can be reinterpreted in the framework of our methodology by making this
choice for $\Phi$, even though the philosophy and derivation in that paper
appear to be very different.  In \citep{Mofjeld2007}, 
the parameters $\xi_0$  and $\sigma$ vary not only
with the tsunami event being considered, but also with the specific point
$(x,y)$ in the region of interest.    In the $\Delta t$ Method or the  Pattern Methods 
we recommend,
the $\Phi$ might be chosen differently for different tsunami events 
(adjusting $\Delta t$ or the pattern), but the same $\Phi(\hat\xi)$ is used for
all $(x,y)$.  

We now summarize the description of our general methodology.  The
key ingredients are the functions $Z$ (and hence the inverse $Z^{-1}$) 
and $\Phi$.  Once these have been chosen, we estimate the probability that
the QoI will exceed some level $\hat\zeta$ by
\eql{Pzeta}
P[\zeta > \hat\zeta] \approx \Psi(\hat\zeta) \equiv 
\Phi\left(Z^{-1}(\hat\zeta)\right).
\end{equation} 
In practice this is applied for specific values of $\hat\zeta$, namely the
exceedance values $\zeta_i$ of interest, in order to compute 
$\Psi(\zeta_i) = \Phi(\hat\xi_i)$ where $\hat\xi_i = Z^{-1}(\zeta_i)$.
In words, this means that for each exceedance level $\zeta_i$ we estimate the
static tide level $\hat\xi_i$ above which $\zeta >
\zeta_i$ in a tsunami simulation, and then we evaluate $\Phi(\hat\xi_i)$ to
determine the probability that the tide will be sufficiently high when the
tsunami ``arrives'' (in the sense as encapsulated in the choice of 
$\Phi$).

\begin{figure}[t]
\hfil\includegraphics[width=4.2in]{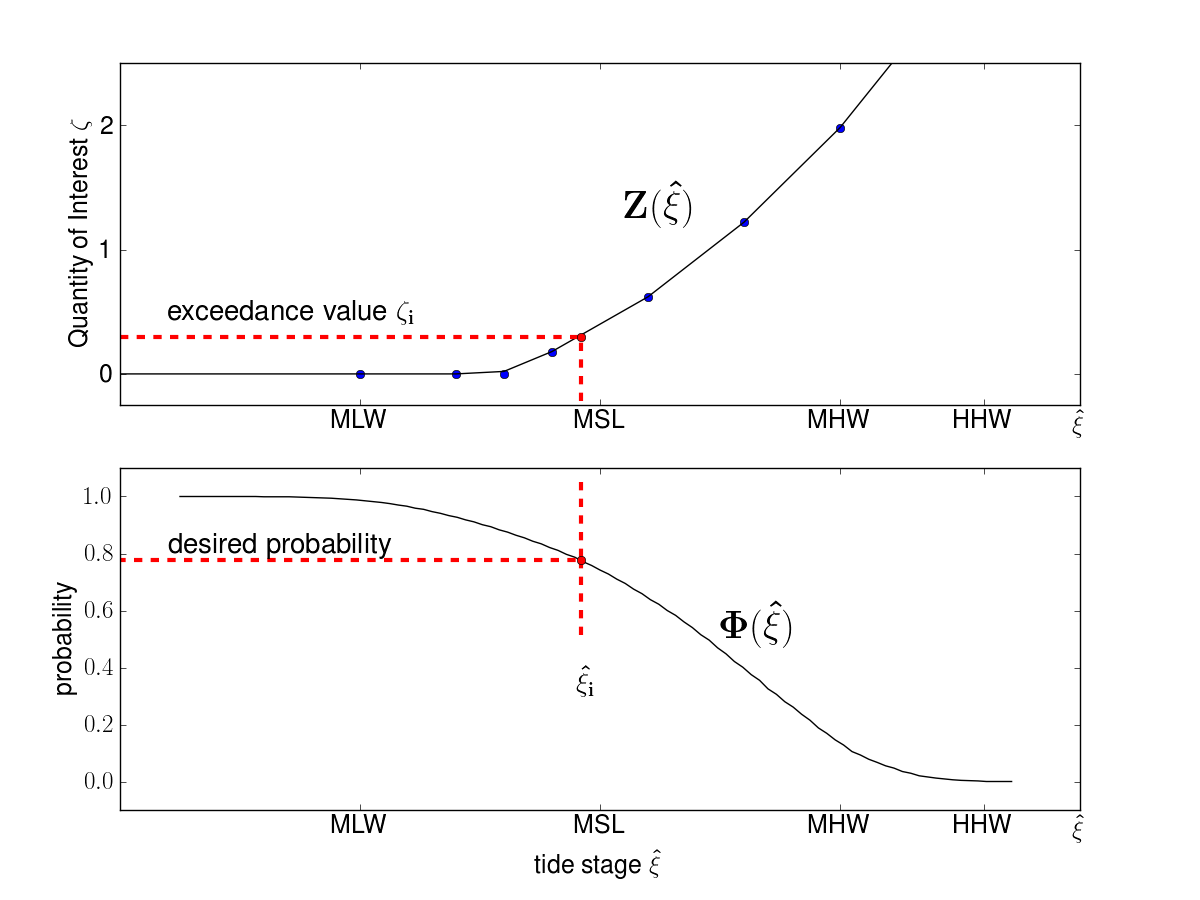}\hfil
\vskip -4mm
\caption{\label{fig:PsiFromZandPhi} Determining $\Psi(\zeta_i)$ 
by computing $\hat\xi_i = Z^{-1}(\zeta_i))$
and evaluating $\Phi(\hat\xi_i)$. }
\end{figure}

This is further illustrated  graphically in \Fig{PsiFromZandPhi}.  
The top figure shows the function $Z(\hat\xi)$.  For a typical $\zeta_i$ on
the vertical axis we invert $Z$ to find $\hat\xi_i$ on the horizontal axis.
Then the bottom figure shows how this same value $\hat\xi_i$ is used to
determine the desired probability by use of the function $\Phi$.

Finally, we note that if $Z(\hat\xi)$ is not a monotone function, then rather
than using a single value $Z^{-1}(\hat\zeta)$ defined by \eqn{Zinv}, we could 
instead determine the intervals $(\hat\xi_j, \hat\xi_{j+1})$ over which
$Z(\hat\xi) > \hat\zeta$ and sum $\Phi(\hat\xi_j)
- \Phi(\hat\xi_{j+1})$ over all such intervals to define $\Psi(\hat\zeta)$.

\section{Methods based on the framework} \label{sec:methods}
We need to estimate $P[\zeta > \hat\zeta]$ in equation (\ref{Pzeta}), the probability that the
QoI $\zeta$ exceeds $\hat\zeta$ whenever a given tsunami event occurs.
We select a method for doing this by choosing a function $Z$ and a function $\Phi$,
such as one of those depicted in Figure \ref{fig:PhiDt}.
The choices of these functions are given in Section \ref{sec:dtmethod} for the 
$\Delta t$ Method, in Section \ref{sec:tides_pattern} for the Pattern Method and in
Section \ref{sec:Gmethod} for the G Method.  These methods are then compared in
Section \ref{sec:tides_compare} using tsunamis from the PTHA study in \citep{GonLevAda}.

\subsection{The $\Delta t$ Method}\label{sec:dtmethod}

This method has been described briefly in Section \ref{sec:overview}. 
The function $Z(\hat\xi)$ is chosen using the second approach outlined in Section \ref{sec:overview}.
We typically use at least three values of $\hat\xi$ including $\xiMLLW$, $\xiMSL$, and $\xiMHHW$
as sealevel
parameters for GeoClaw simulations of the shallow water equations.  The resulting $\hat{Z}$
values are used to make the piecewise linear $Z$ function shown in red in Figure \ref{fig:PhiDt}.

Once $\Delta t$ has been selected for a particular tsunami,
the function $\Phi(\hat\xi)$ is given in equation (\ref{PhiDt}) and its graph
can be found in Figure \ref{fig:PhiDt}.
Discrete values of $\Phi(\hat\xi)$ are gotten by first placing valid values of $\hat\xi$ in bins and
then using a $\Delta t$-window of time and sliding it one minute at a time across a year's worth
of tide gauge data at the tsunami destination site of interest.
Each time the $\Delta t$-slider window stops, we find the maximum
tide level within the window. We increment a counter in the first bin whose right edge
exceeds or equals this maximum (to create a histogram) and also in all lower bins
(to create a cumulative histogram).  Dividing by the number of times the
$\Delta t$-slider window stops gives the probability mass function and $\Phi(\hat\xi)$,
respectively. The PDF $\phi(\hat\xi)$ results by dividing
the probability mass function by the bin size.
The $\Phi$ values for each bin's left edge are stored in a table and interpolated as needed.

Different
tsunamis will require different choices of $\Delta t$.
We place
computational gauges at various locations
where we collect time series output to determine the width and relative occurrence times of
potentially damaging waves.
The width of the responsible wave of
biggest amplitude certainly gives a minimum value for the {\em contiguous} $\Delta t$ interval,
and we increase $\Delta t$ if there are nearby waves of nearly equal amplitude, so the tsunami is
effectively modelled as one square wave of width $\Delta t$.

In Section \ref{sec:tides_compare}, we see the $\Delta t$ Method works remarkably
well compared to the Pattern Method for appropriately chosen $\Delta t$, and there
we give recommended values for particular tsunamis in the Crescent City study.

\noindent
\subsection{The Pattern Method}\label{sec:tides_pattern}

$Z(\hat\xi)$ is chosen 
the same as for the $\Delta t$ Method, see Section \ref{sec:dtmethod}.
When the tsunami consists of only one wave, we will see that the
Pattern Method is simply the $\Delta t$ Method where $\Delta t$ is
the time duration of the wave.

The Pattern Method uses the relative heights of the
waves seen at a computational gauge located in the water,
their widths, and the relative times they occurred,
to create the $\Phi(\hat\xi)$
associated with this particular wave pattern. 
This is extra work, but the difference
is that a {fixed} $\Delta t$ will not have to be chosen.
Instead, the entire pattern will be taken into account to calculate $\Phi(\hat\xi)$.

Suppose the tsunami has $K$ waves.  We model wave $W_k$ with a
square wave and record the difference of its height from that of
the highest wave as $D_k$. That is, $D_k = H - H_k$ where $H_k$ is the height
of wave $W_k$, and $H=\max_k H_k$.
We record the starting and terminating
times of $W_k$ as the interval $I_k = [S_k , T_k ]$.  These times
are relative to the start of $W_1$, so we set $S_1 = 0$; they are
recorded in minutes since our tide record
has minute data. The pattern's duration is then
$T_K$ minutes.

As an example, in Figure \ref{fig:AASZe02gauge101}, we show the GeoClaw tsunami for an
Alaskan Aleutian event  
(AASZe02) in red and the pattern in black.
The first wave arrived at Crescent City 4 hours and 23
minutes after
the earthquake and nothing significant was seen there after 11 hours.  The pattern
is well represented by the 7 waves shown. We are overestimating the probability
a bit by using square waves, but we don't have to account for tides between these waves.
Table \ref{tab:patternvalues} shows the values that describe the pattern. We note that
the first wave began at 263 minutes after the earthquake and the amplitude of the
largest wave $W_7$ was about 1.5 meters. The black horizontal line starts at 0.2 meters
since the GeoClaw run was done at $\xiMHHW$ which is 0.2 meters above $\xiMHW$, the zero level
for the plot in Figure \ref{fig:AASZe02gauge101}.

\vskip 2mm
\begin{tabular}{ll}
\begin{minipage}{2.0in}
\includegraphics[width=2.0in]{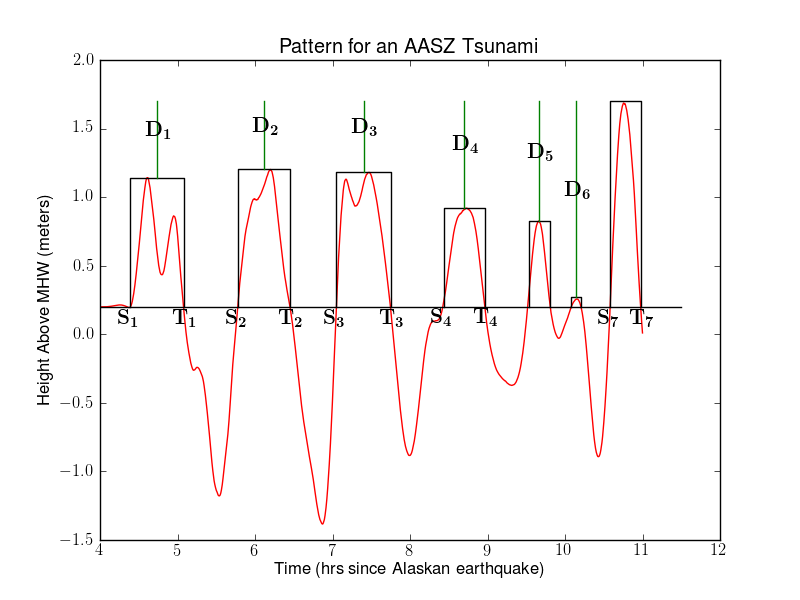}
\vskip -3mm
\captionof{figure}{\label{fig:AASZe02gauge101} Pattern for AASZe02}
\end{minipage} &
\begin{minipage}{2.0in}
\captionof{table}{{Pattern Values}
\label{tab:patternvalues}}
\vskip -1.3mm
\small\begin{tabular}{||c|c|c||} \hline
  Wave                  & $I_k=[S_k, T_k ]$         & $D_k$ (meters)   \\ 
  $W_k$                 & Wave Interval            & Difference to  \\ 
                        & (min since $S_1$)       & Tallest Wave  \\ \hline
  $W_1$                 & [000, 042]     & 0.561          \\ \hline 
  $W_2$                 & [084, 124]     & 0.498           \\ \hline 
  $W_3$                 & [160, 202]     & 0.517           \\ \hline 
  $W_4$                 & [243, 275]     & 0.782           \\ \hline 
  $W_5$                 & [309, 325]     & 0.876           \\ \hline 
  $W_6$                 & [342, 349]     & 1.450            \\ \hline 
  $W_7$                 & [372, 396]     & 0.000            \\ \hline 
\end{tabular}
\end{minipage}\\
\end{tabular}

\vskip 2mm
As in the $\Delta t$ Method, we 
put the valid values of $\hat\xi$ into bins.
But now we take our pattern-slider window that has length $T_K$ and
slide it one minute at a time across a year's worth
of the tidal record.  Each time the window stops we
calculate $\max_{1 \leq k \leq K} ( \,\max_{t \in I_k} \xi(t + t_0 ) - D_k \, )$ corresponding
to a different tsunami start time $t_0$. Then we
increment a counter in the first bin whose right edge exceeds or equals this value (and
in all lower bins)
to create a histogram (cumulative histogram). These histograms are used as before to
determine $\Phi$ and $\phi$ values.
Thus, the Pattern Method $\Phi$
function is
\begin{equation}\label{PatPHI}
\Phi_{\rm Pattern} (\hat\xi) = P[ \max_{1 \leq k \leq K} ( \,\max_{t \in I_k} \xi(t + t_0 ) - D_k \, ) > \hat\xi ].
\end{equation}

Note, that if the pattern consists of only one wave, then $D_1 = 0$
and the Pattern Method is just the $\Delta t$ Method with $\Delta t$ being
the length of $I_1$.
The Pattern Method has advantages over the $\Delta t$ Method.  
Only one synthetic gauge needs to be examined.
A wave with amplitude less than the
maximum one could also cause exceedance of $\hat\zeta$ if it
occurred
at a time when the tide level was sufficiently high.
Tsunamis with longer duration are more accurately represented
since the tide record during each interval $I_k$ and not between
needs to be examined. This gives an automatic procedure that avoids
the difficulty in choosing an appropriate $\Delta t$.
A possible limitation is that
the Pattern-Method requires the simulation code to have GeoClaw's capability of a computational gauge.

\subsection{The G Method}\label{sec:Gmethod}
We describe how this method fits our framework in Section \ref{sec:overview}.  Only
one GeoClaw simulation at $\xiMHHW$ is done to get $\hat{Z}(\xiMHHW)$.
The function $Z(\hat\xi)$ for location $(x,y)$ within the destination area
of interest
is then given in equation (\ref{Zlinear}).

A proxy tsunami is assumed for each $(x,y)$, defined as having a duration of 5 days (or
$T_G = 7200$ minutes)
with e-folding
time of 2 days and period of 20 minutes with maximum amplitude $A_G = A_G (x,y)$.
This assumed tsunami has maximum height $H$ which
can be measured to any
fixed reference level for our purposes and occurs with the first wave having amplitude
$A_G$.  Its height at time $t$ after the first wave begins is $H(t)$, and following
previous notation, the distance to the maximum is $D(t)=H-H(t)$. 
Then $\Phi_G (\hat\xi)$ is 

\begin{equation} \label{PhiGnoapprox}
\Phi_G (\hat\xi) = P \, [\, \max_{t_0 \leq t \leq t_0 + T_G} (\xi(t) - D(t))\,] > \hat\xi ,
\end{equation}

\noindent
and can be further approximated by equation (\ref{PhiErf}) with mean $\xi_0$ and
standard deviation $\sigma$. Using $\hat\zeta = Z(\hat\xi)$ and (\ref{Zlinear}), we get
$\hat\xi = \hat\zeta - \hat{Z}(\xiMHHW) + \xiMHHW$ and
(\ref{PhiErf}) becomes
\begin{equation}\label{zetaerf}
P[\zeta > \hat\zeta] \approx \frac{1}{2} \left ( 1- \rm{erf}\left ( (\hat\zeta - \zeta_0 )/\sqrt{2}\sigma \right ) \right )
\end{equation}
\noindent
where the mean $\zeta_0$ of $\zeta$ is given by
\begin{equation} \label{zeta0}
\zeta_0 = \hat{Z}(\xiMHHW) - \xiMHHW + \xi_0 .
\end{equation}
\noindent
This makes sense. As the value $\xi_0$ approaches $\xiMHHW$, the mean of $\zeta$ should be
$\hat{Z}(\xiMHHW)$. Also, $\zeta$'s value is $\hat{Z}(\xiMHHW)-\xiMHHW$
when $\hat\xi = 0$, or $\xiMSL$.  As $\xi_0$ approaches 0, $\zeta_0$ should
approach this value.

We can now make the final connection to the formula for $P[\zeta > \hat\zeta]$ given in
\citep{Mofjeld2007}.  There, $\xi_0$ and $\sigma$ are functions of location $(x,y)$ and
are given by
$\xi_0 = C \xiMHHW e^{-\alpha {(A_G / \sigma_0 )}^\beta}$,
and $\sigma = \sigma_0 (1 - C^{\prime} e^{-\alpha^\prime {(A_G / \sigma_0 )}^{\beta^\prime}})$.
When $\zeta + \xi_{ref}$ is the flow depth above $\xiMSL$, and $S$ is the amount of subsidence or uplift
(positive with subsidence so that $\xiMHHW - S$ is the subsided background water of the
simulation),
it makes sense to express $A_G$ as
\begin{equation}\label{amplitude}
A_G = \hat{Z}(\xiMHHW) +\xi_{ref} - (\xiMHHW - S).
\end{equation}

\noindent
Substituting into (\ref{zeta0}) gives
\begin{equation} \label{zeta0Aform}
\zeta_0  = A_G  -\xi_{ref} - S + C \xiMHHW e^{-\alpha {(A_G /\sigma_0)}^\beta}.
\end{equation}

\noindent
Using $\xi_{ref}=\xiMLLW$ and $S=0$ gives the method in \citep{Mofjeld2007},
where the
parameters are also given for a variety of tsunami destinations.
Those for Crescent City include $\sigma_0 = 0.638$, the standard deviation for the
tides there, and the regression parameters
$\alpha^\prime$ = 0.056,
$\beta^\prime$ = 1.119,
$C^\prime$=0.707,
$\alpha$=0.17,
$\beta$=0.858, and
$C$=1.044.

The G Method has two major limitations. First,
only one GeoClaw simulation with tide level $\hat{\xi} = \xiMHHW$ is used.
This is appropriate whenever $Z(\hat\xi)$ is indeed a linear function of
slope 1, since then $A_G$ will be uniquely defined from only one simulation.
The $Z(\hat\xi)$ functions created with multiple simulations show this is not true,
especially 
at onshore locations, see Figure \ref{fig:DepthLocations}.

The second limitation is the use of the same 5-day proxy tsunami (where the amplitude alone varies
at each $(x,y)$)
as the pattern for modelling each tsunami in a PTHA study, especially when the major question
being studied is the flow depth at land points.  As seen in Table \ref{MofPattern}, the duration of all the
tsunamis studied that could impact the maximum inundation at a land point is much less than 5 days,
and we observe these tsunamis have patterns that are very different.
A local tsunami from the Cascadia Subduction Zone will typically have only one or two waves occuring
over a short time frame that are responsible for the maximum; whereas, far field events can have
damaging waves occuring over a longer time frame whose amplitudes can
{\em increase} during significant tidal variations.
\section{An example PTHA study and method comparisons} \label{sec:tides_compare}

\subsection{The Crescent City PTHA study}\label{sec:CCtides}

Phase I of the PTHA study of Crescent City, California reported in \citep{GonLevAda}
focused on flow depths. Output products such as 100- and
500-year hazard maps, hazard curves at specific locations, and probability contours for exceeding
a specific $\hat\zeta$ level can be found in the report. All these products used the
Pattern Method. Here we demonstrate why this is the preferred method for including tidal
uncertainty.

GeoClaw simulations of the shallow water equations were conducted at multiple static tide
levels for each tsunami in the study to find the QoI at each fixed grid location.  For this
study, the QoI was the maximum flow depth above topography/bathymetry (onshore points) and the maximum
flow depth plus the original bathymetry (offshore points) measured in meters and was denoted $\zeta$.
For offshore points the bathymetry is negative, and represents the negative of the distance
between the underwater topography and $\xiMHW$. The QoI for offshore points is then
the amount of flow
depth above $\xiMHW$ plus the amount of subsidence (original minus final bathymetry).
The QoI for onshore points is the flow depth measured above the final topography. With this
definition, the QoI is continuous at the shoreline which corresponds to $\xiMHW$.

The GeoClaw simulations made use of computational gauges placed at strategic locations where
the time series of the QoI were monitored.  One such gauge, called Gauge 101, was placed in
the Crescent City harbor and was where some of our results are reported. In particular, this
gauge is where we recorded each tsunami's pattern for the Pattern Method.

\begin{figure}[h]
\vskip -2mm
\hfil\includegraphics[width=1.5in]{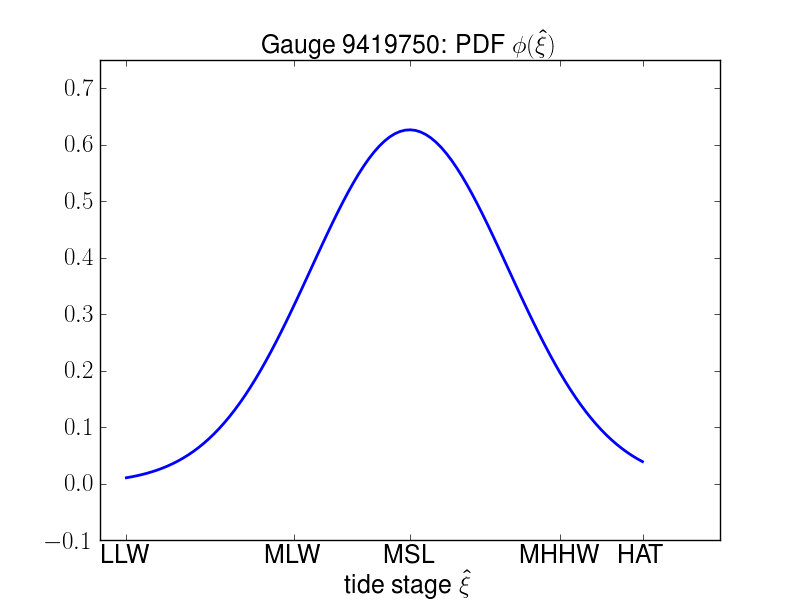}\hfil
\hfil\includegraphics[width=1.5in]{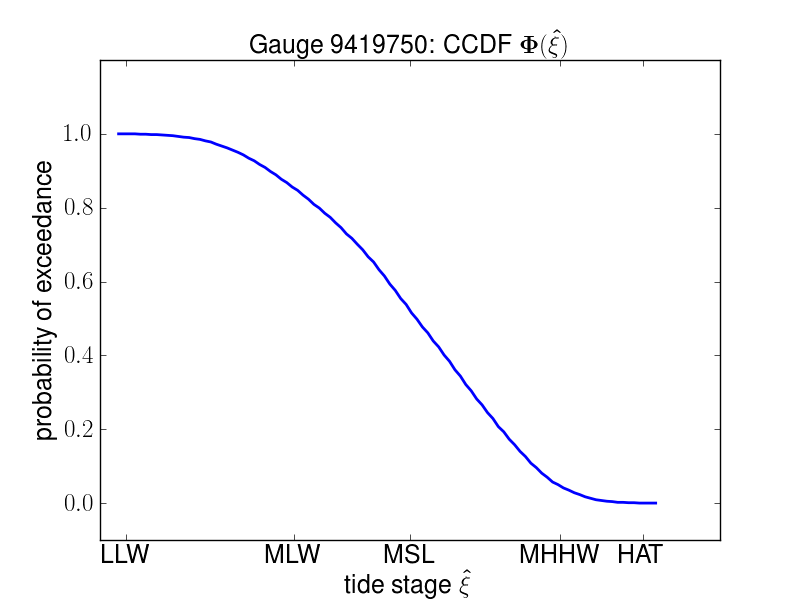}\hfil
\vskip -1mm
\caption{\label{fig:CCpdfandcumul} Crescent City Tides: 
Left: $\phi_0 (\hat\xi)$ ($\xi_0=0.0$, $\sigma_0$ = 0.638),
Right: $\Phi_0 (\hat\xi)$}
\end{figure}

A few important tidal constants at Crescent City Gauge No. 9419750, see
\url{http://tidesandcurrents.noaa.gov/waterlevels.html?id=9419750},
that were used as sealevel parameters for GeoClaw
simulations are 
Mean Lower Low Water ($\xiMLLW=-1.13$), Mean Low Water ($\xiMLW=-0.75$),
Mean Sea Level ($\xiMSL=0.0$), Mean High Water ($\xiMHW=0.77$), and Mean Higher High 
Water ($\xiMHHW=0.97$).
The lowest and highest
water seen at the gauge in a year's data from July 2011 to July 2012 are
$\xi_{Lowest}=-1.83$ and $\xi_{Highest}=1.50$ in meters, referenced to $\xiMSL$ respectively. 
Figure \ref{fig:CCpdfandcumul} shows the PDF $\phi_0 (\hat\xi)$ and the CCDF $\Phi_0 (\hat\xi)$
for this yearly data.  $\Phi_0 (\hat\xi)$ is that of the $\Delta t =0$ Method.

Column 1 of Table \ref{MofPattern} gives representative tsunamis used in this PTHA study.
The acronyms CSZ, AASZ, KmSZ, KrSZ, and SchSZ stand for the Cascadia, Alaskan Aleutian, Kamchatka,
Kuril, and South Chile Subduction Zones, respectively, and TOH refers to Tohoku.  We also denote
tsunami events in the form AASZe03; for example, event number 3 on the Alaska Aleutian Subduction Zone.
Some events, e.g., a CSZ Mw 9.1 event, have multiple realizations.
CSZBe01r01-CSZBe01r15 refers to the CSZ Bandon sources modeled as 15 realizations of different
slip distributions for a single event used in a PTHA study of Bandon, Oregon
\citep{Witter:BandonSP43}. More details about these earthquake source models can be found in
\citep{GonLevAda}.

The recommended value of $\Delta t$ for all tsunamis in Table \ref{MofPattern} and in \citep{GonLevAda}
can be given.  We recommend $\Delta t=1$ for the
Kamchatka event  KmSZe01 and $\Delta t=3$ for KmSZe02.  For the three
Kuril events, we recommend
$\Delta t=2$ for KrSZe01, $\Delta t=3$ for KrSZe02, and $\Delta t=4$ for KrSZe03. For the Alaska events,
we recommend $\Delta t=1$ with the exception of $\Delta t=2$ for AASZe02.  The value $\Delta t=1$ should
be used for the Chilean event SChSZ01, the Tohoku event TOHe01,
and the Cascadia Bandon CSZBe01r13 and CSZBe01r14 tsunamis.  The value $\Delta t=0$
should be used for the remaining Cascadia Bandon tsunamis, CSZBe01r01-CSZBe01r12
and CSZBe01r15.

For the $\Delta t$ and Pattern Methods most of these tsunami events were simulated by running
GeoClaw at $\xiMLLW$,
$\xiMSL$, and $\xiMHHW$, respectively.
For more intense analysis, the AASZe03 event (similar
to the 1964 Alaska tsunami) was run using 11 tide levels.
These levels referenced to $\xiMSL$ were  -1.13, -0.75, -0.50, -0.25, 0.00,
0.25, 0.50, 0.77, 0.97, 1.25, and 1.5 meters.
For the G Method, only the GeoClaw $\xiMHHW$ results were required.

\begin{figure}[t]
\hfil\includegraphics[width=2.3in]{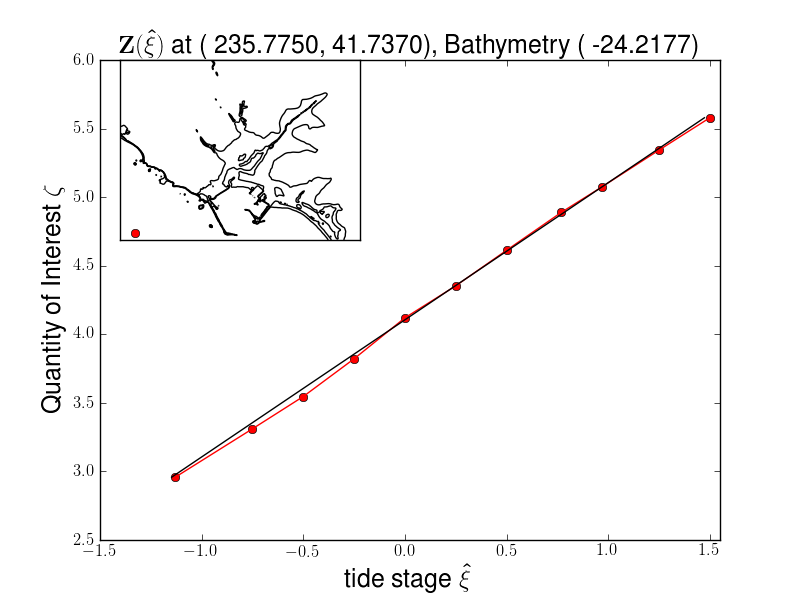}\hfil
\hfil\includegraphics[width=2.3in]{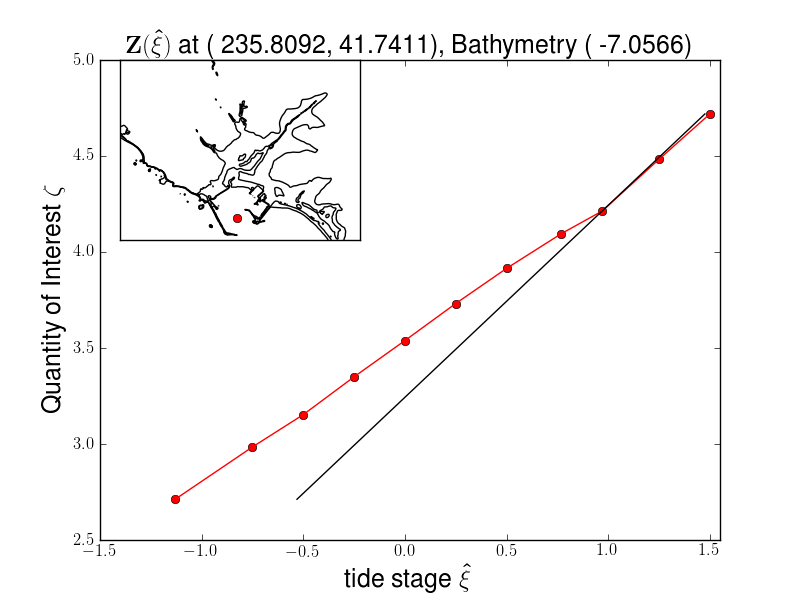}\hfil
\hfil\includegraphics[width=2.3in]{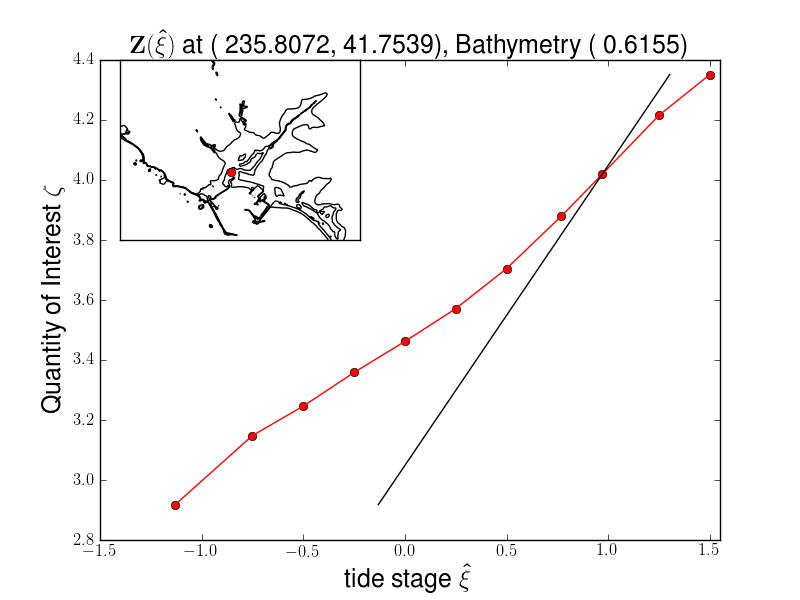}\hfil
\hfil\includegraphics[width=2.3in]{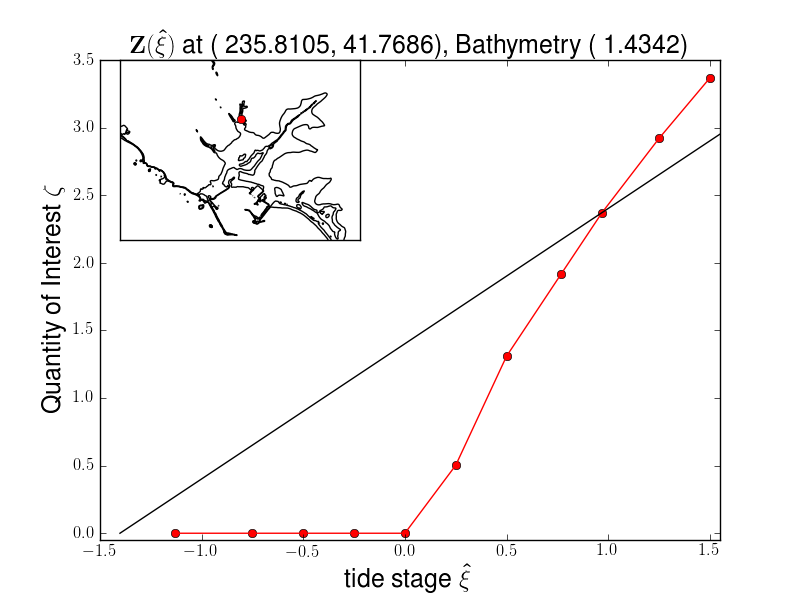}\hfil
\vskip -3mm
\caption{\label{fig:DepthLocations} $Z (\hat\xi)$ functions for flow depth
at land locations (bottom row) and flow depth plus bathymetry at offshore locations (top row).
Location longitude, latitude, and bathymetry are given for each location.}
\end{figure}

Figure \ref{fig:DepthLocations} gives the $Z(\hat\xi)$
functions for four different locations for the AASZe03 tsunami.
The black line on the plots is the $Z(\hat\xi)$
used with the G Method and is the slope 1 line through the point corresponding to
$\hat\xi=\xiMHHW=0.97$. The red line is the $Z(\hat\xi)$ used by the $\Delta t$ and Pattern Methods.
The longitude, latitude, and
bathymetry of the location is given on the graphs.
The top row shows
the $Z(\hat\xi)$ function for two offshore locations is similar for all three methods.
The second row
shows the $Z(\hat\xi)$ function for the $\Delta t$ and Pattern Methods
for two onshore locations
is far from the slope 1 line used by the G Method. 

In Section \ref{multiple}, we compare the Pattern and G Methods based on the mean and standard
deviations of their $\phi$ PDFs for many tsunamis. Graphs of the $\Phi$ CCDFs for the AASZe02
and AASZe03 tsunamis are given in Section \ref{dtPat} for the $\Delta t$ and Pattern Methods. 
In Section \ref{validation}, we use graphs of the $\Phi$ functions for the Pattern and G Methods
for the AASZe03 tsunami and vary the duration of its Proxy (used by the G Method) as a validation
of the Pattern Method.  Finally, in
Section \ref{PTPD}, we use the AASZe03 tsunami to compare
all three methods, including their probabilities for
exceeding the 35
values $\zeta_i$ = 0, 0.1, 0.2, $\dots$, 1.9, 2.0, 2.5, $\dots$, 5.5, 6.0,
7.0, $\dots$, 12.0 when the tsunami occurs.

\subsection{G and Pattern $\phi$ comparisons at Gauge 101 for multiple sources}\label{multiple} 
In Table \ref{MofPattern}, we compare the probability density functions $\phi$ of the G and
Pattern Methods
for some of the tsunamis considered in the Crescent City study.
Table \ref{MofPattern} shows there are significant differences between the G Method
and the Pattern Method.
Only for the five large amplitude tsunamis CSZBe01r01,
CSZBe01r02, CSZBe01r03, CSZBe01r04, and CSBe01r05 do the two methods have $\phi$'s with similar means
and standard deviations.  For the other tsunamis in the table, the G Method has a
much higher mean and smaller standard deviation than the Pattern Method.  

\begin{table}[h]
\begin{center}
\caption{G and Pattern Method $\phi$ comparisons at the Gauge 101 location.
The length $T$ in min. and amplitudes $A_G=\hat{Z}(\xiMHHW)+\xiMHW-\xiMHHW$
in m. are given in columns 2 and 3 for some tsunamis used
in this study.  Columns 4-7 give the mean $\xi_0$ and standard deviation $\sigma$
for these methods.
}

\label{MofPattern}
\begin{tabular}{||l|r|r|r|r|r|l||} \hline
           &      &       &G            &G        & Pattern   &Pattern \\
Source     &T     &$A_G$  &$\xi_0$      &$\sigma$ & $\xi_0$     &$\sigma$ \\
Name       &(min) &(m)    &(m)       & (m)     & (m)       &(m)      \\ 
           &      &       &          &         &           &          \\ \hline
AASZe03-Proxy   &7205  &3.92   &0.45       &0.34      & 0.46       &0.34       \\ \hline
AASZe01	   & 328  &1.96   &0.65       &0.27      & 0.12       &0.53       \\ \hline
AASZe02	   & 396  &1.50   &0.71       &0.25      & 0.36       &0.37       \\ \hline
AASZe03	   & 267  &3.92   &0.45       &0.34      & 0.14       &0.60       \\ \hline
%
%
AASZe08	   & 114  &0.30    &0.93       &0.20      &0.18       &0.60       \\ \hline
KmSZe01	   & 308  &0.92    &0.80       &0.22      & 0.15       &0.54       \\ \hline
KrSZe01	   & 275  &0.50    &0.88       &0.21      &0.22       &0.52       \\ \hline
%
%
SChSZe01   & 106  &0.60    &0.86       &0.21      &0.16       &0.60       \\ \hline
TOHe01     & 324  &1.66   &0.69       &0.26      &0.07       &0.59       \\ \hline
CSZBe01r01 & 329  &14.18  &0.09       &0.56      &0.04       &0.63       \\ \hline
CSZBe01r02 & 326  &12.96  &0.11       &0.55      &0.04       &0.63       \\ \hline
CSZBe01r03 & 326  &13.31  &0.10       &0.55      &0.04       &0.63       \\ \hline
CSZBe01r04 & 157  &13.00  &0.11       &0.55        &0.04       &0.63       \\ \hline
CSZBe01r05 & 160  &11.30  &0.14       &0.53      &0.04       &0.63       \\ \hline
CSZBe01r07 & 160  &7.78   &0.24       &0.46      &0.03       &0.63       \\ \hline
CSZBe01r08 & 161  &6.56   &0.29       &0.43      &0.03       &0.63       \\ \hline
%
CSZBe01r10 & 160  &2.39   &0.60       &0.29      &0.03       &0.63       \\ \hline
CSZBe01r11 & 163  &4.79   &0.39       &0.37      &0.03       &0.63       \\ \hline
%
\end{tabular}
\end{center}
\end{table}

We note that the $\Delta t$ method with a good choice of $\Delta t$ gives
very similar results to the Pattern Method and was not included in Table
\ref{MofPattern}. For example, for the AASZe03 event, the
Pattern Method values were $\xi_0 = 0.14$ and $\sigma =0.60$. These same values
for the $\Delta t$ method were 0.12 and 0.62, respectively.

\subsection{$\Delta t$ and Pattern $\Phi$ comparisons for AASZe03 and AASZe02}\label{dtPat}
\begin{figure}[h]
\hfil\includegraphics[width=2.0in]{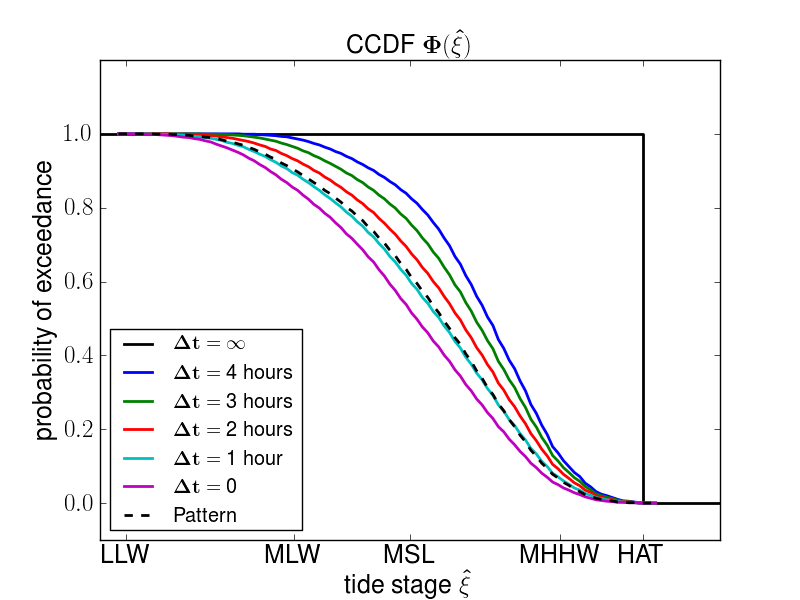}\hfil
\hfil\includegraphics[width=2.0in]{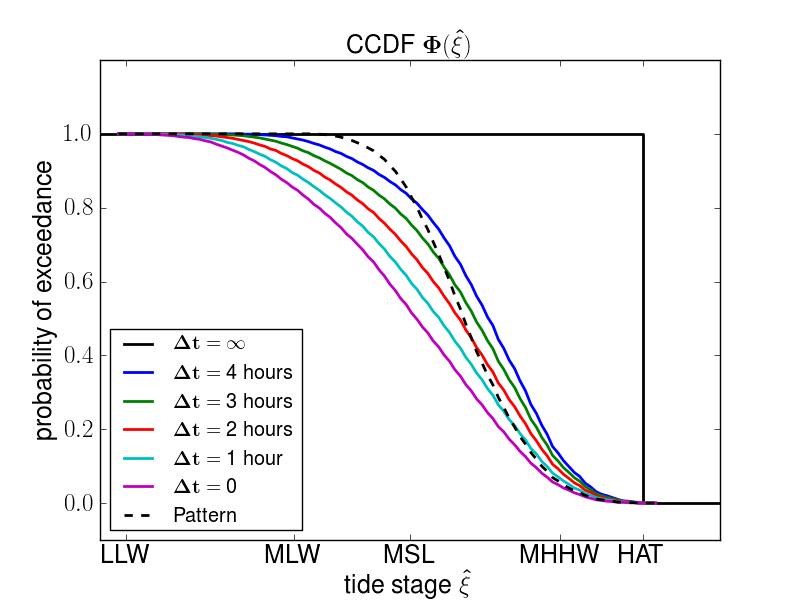}\hfil
\vskip -3mm
\caption{\label{figAASZcumul1} $\Phi(\hat\xi)$ Comparisons.
Left: AASZe03, Right: AASZe02
}
\end{figure}

Figure \ref{figAASZcumul1} shows that for some tsunamis $\Phi_{\rm{Pattern}}(\hat\xi)$
is similar to
$\Phi_{\Delta t}(\hat\xi)$
for a fixed $\Delta t$ (AASZe03, $\Delta t=1$), while for other tsunamis $\Phi_{\rm{Pattern}}(\hat\xi)$
is consistent with a varying $\Delta t$ (AASZe02). 
$\Phi_{\rm{Pattern}}$ is shown as
a dotted line on the same graph as the $\Phi_{\Delta t}$'s for varying $\Delta t$.

\subsection{G and Pattern $\Phi$ comparisons at Gauge 101 for AASZe03}\label{validation}
\begin{figure}[h]
\hfil\includegraphics[width=2.2in]{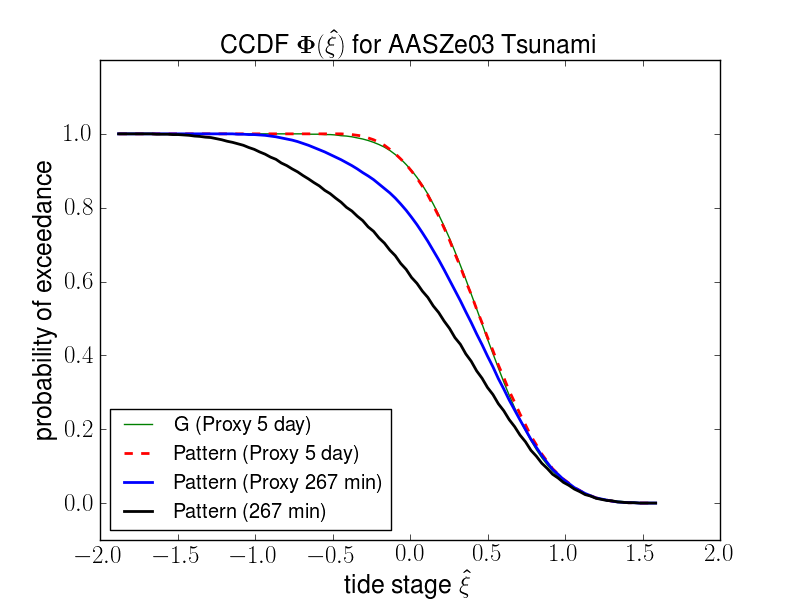}\hfil
\vskip -2mm
\caption{\label{fig:validate_pattern} Pattern Method Validation}
\end{figure}

We ran the Pattern Method on the 5-day proxy tsunami that is assumed by the G Method and compared
the resulting $\Phi$ functions at Gauge 101. The amplitude for the 5-day proxy
tsunami 
was taken as that of the biggest wave seen at Gauge 101 for AASZe03.
The two distributions when plotted are almost
identical with values differing mostly less than 1\% as seen in Figure \ref{fig:validate_pattern}
as the green and dashed red graphs and given in the first line of Table \ref{MofPattern}.
The black graph is the distribution for the
Pattern Method on the actual tsumani at Gauge 101 for which we used a $T=267$ minute duration.
The blue graph shows the Pattern Method assuming a 267 minute proxy tsunami which shows differences
in the proxy and actual tsunami patterns.
This explains that differences in $\Phi$
for any real tsunami are
not due to our methodology, but to the fact that the real tsunami is not well approximated
by the proxy one, even if we enforce both to have the same time duration.

\subsection{Probability differences} \label{PTPD}
For each grid location, we compared the 35 probabilities
$P( \zeta > \zeta_i )$ of the three methods for AASZe03.
The numbers in Table \ref{tabLoy2} are over all the grid locations that
cover the Crescent City area.
The row labelled {\bf max}
is the maximum difference seen when the method being compared to the Pattern Method
gives the larger result, and the row labelled {\bf min} is the difference
seen when the Pattern Method gives the larger result.
Indeed, differences close to 1 are observed in the first column and
the second column shows that the $\Delta t$ Method (with $\Delta t = 1$) gives
results very close to those of the Pattern Method.
\begin{table}[h]
\begin{center}
\caption{Probability Differences}
\label{tabLoy2}
\vskip -1mm
\begin{tabular}{||r|r|r|} \hline
               &{\bf G - Pattern}    & {\bf $\Delta t$ - Pattern}  \\ \hline 
{\bf max}      & +0.747             & +0.006             \\ \hline 
{\bf min}      & -0.936             & -0.017            \\ \hline
\end{tabular}
\end{center}
\vskip -3mm
\end{table}
Both the $\Delta t$ and Pattern Methods use the amplitude of the
tsunami at Gauge 101 and assume its duration is
T minutes instead of 5 days. Further analysis given in \citep{AdaLevGonz} shows that almost
all of the -0.936 is due to the G Method's choice of 5 days, while all but 0.158 of the 0.747
is due to this choice.  This remaining 0.158 difference is due to the use of a proxy
decaying e-folding pattern of 2 days for the tsunami, rather than the observed
pattern.

\begin{figure}[h]
\hfil\includegraphics[width=2.2in]{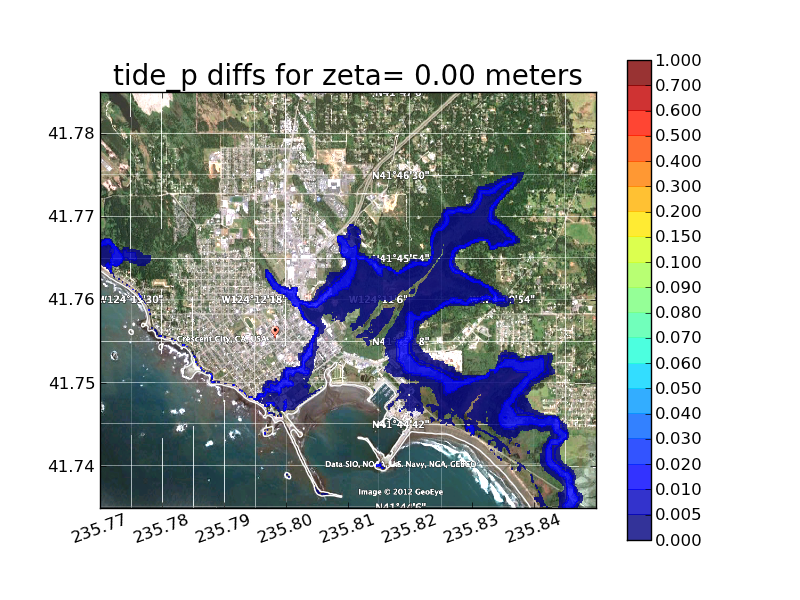}\hfil 
\hfil\includegraphics[width=2.2in]{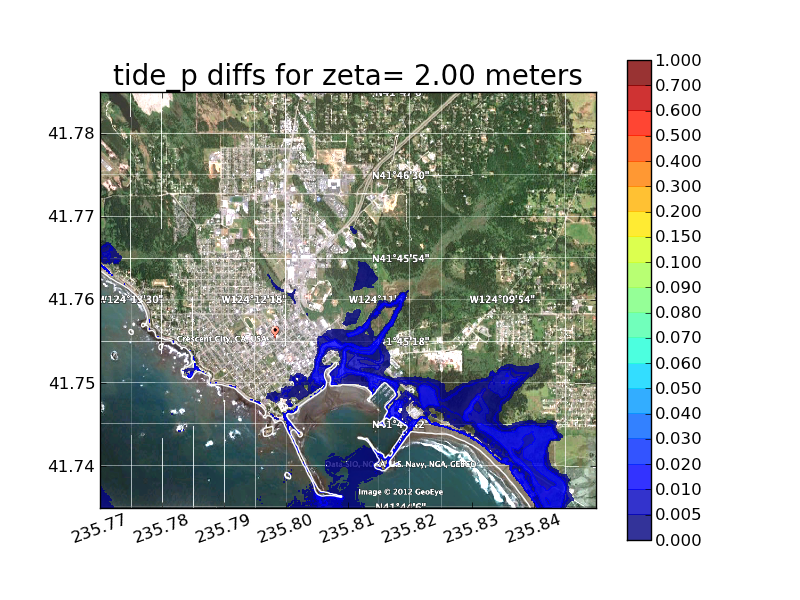}\hfil
\vskip 1mm
\hfil\includegraphics[width=2.2in]{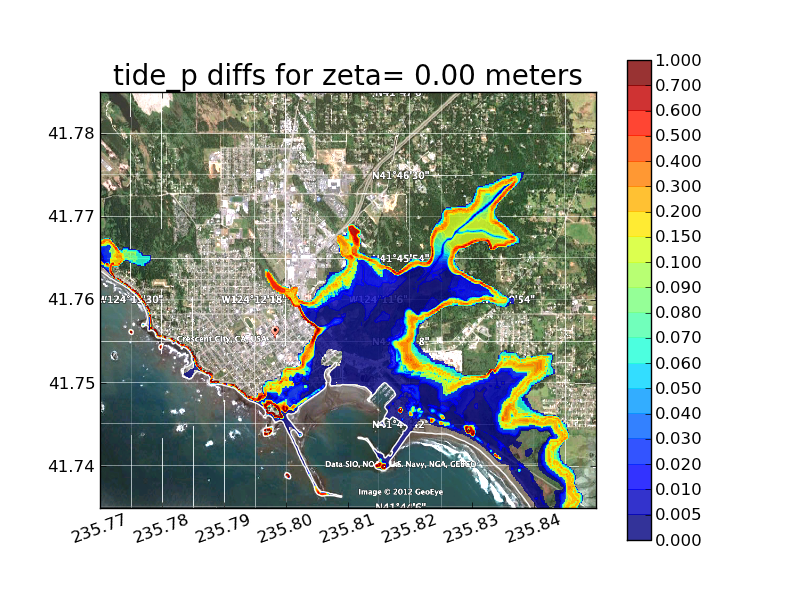}\hfil
\hfil\includegraphics[width=2.2in]{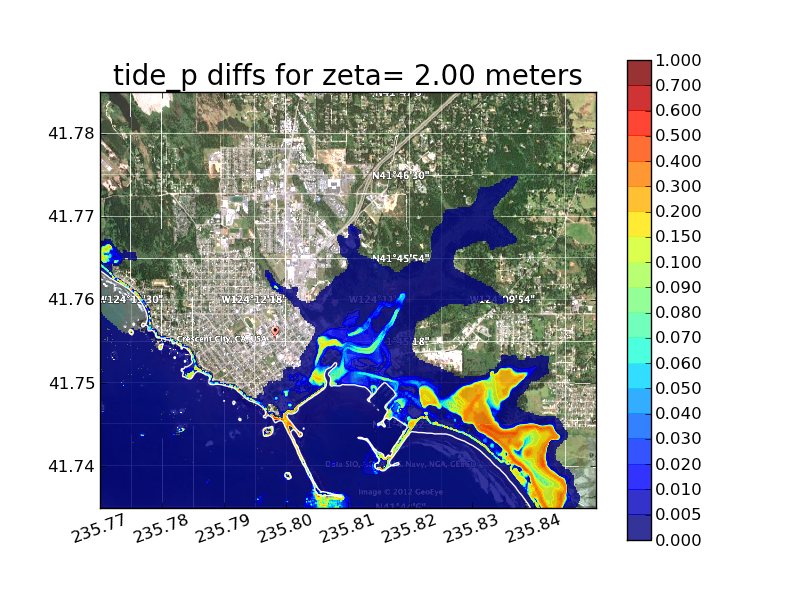}\hfil
\vskip -2mm
\caption{\label{figLoyPDC1} $P(\zeta > \zeta_i )$ Difference Contours, Left: $\zeta_i = 0$ m., Right: $\zeta_i = 2$ m. Top: abs($\Delta t$ - Pattern), Bottom: abs(G - Pattern)
}
\end{figure}

In Figure \ref{figLoyPDC1} we
compare the $\Delta t$ Method and the G Method to the Pattern Method by
giving contour plots of the absolute value of the probability differences of
exceeding $\zeta_i =0$ meters and $\zeta_i = 2$ meters. The brighter colors
in the plots indicate where these probabilities differ the most. As expected 
from Table \ref{tabLoy2}, the $\Delta t$ and Pattern Methods differ less than
2\%.
\section{Conclusions and open questions}\label{sec:conclusions}

The $\Delta t$ Method and the Pattern Method give
quite similar results for a properly chosen $\Delta t$ but vary significantly
from the G Method, especially at land points. 
The Pattern Method
is a very robust method coupled to the wave pattern
for each individual
tsunami and gives modelers a single method that can be used for both
land and water locations.
Both these methods were designed to use GeoClaw simulation
information
at multiple {\em but static} tidal levels, and will work with other
codes that have the capability to produce similar results.

We do not model the currents that are generated by the tide rising and falling.  
A tsunami wave arriving on top of an incoming tide could potentially
inundate further than the same amplitude wave moving against the tidal
current, even if the tide stage is the same.  Modeling this is beyond the
scope of current tsunami models.  

\end{document}